\documentclass[12pt, reqno]{amsart}
\usepackage{hyperref}

\usepackage{amssymb,amsmath,graphicx,amsthm}

\def\contrazione{\raisebox{1pt}{\,{\mbox{\tiny{$|\!\raisebox{-0.7pt}{\underline{\hphantom{X}}}$}}}\,}}
\def\Op{\T{Op}^{\T{ord}(\Psi)-\frac12}}
\def\Opbis{\T{Op}^{2\T{ord}(\Psi)-1}}
\def\sumKk-1{\underset{|K|=k-1}{{\sum}'}}
\def\sumKq-1{\underset{|K|=q-1}{{\sum}'}}

\def\sumJ(p-1){\underset{|J|=p-1}{{\sum}'}}

\def\sumKq{\underset{|K|=q}{{\sum}'}}
\def\sumKp-2{\underset{|K|=p-2}{{\sum}'}}

\def\sumjq+1{\underset {j\leq q+1}\sum}
\def\sumjn-1{\underset {j\leq n-1}\sum}

\def\sumh{\underset {h=1}{\overset {n-1}\sum}}

\def\epf{ \hskip15cm$\Box$}
\def\bpf{{\it Proof. }\hskip0.2cm}
\def\bt{\begin{theorem}}
\def\el{\end{lemma}}
\def\bl{\begin{lemma}}
\def\et{\end{theorem}}
\def\bp{\begin{proposition}}
\def\ep{\end{proposition}}
\def\bd{\begin{definition}}
\def\ed{\end{definition}}
\def\br{\begin{remark}}
\def\er{\end{remark}}

\def\simleq{\underset\sim<}

\def\simgeq{\underset\sim>}

\def\T{\text}

\def\1#1{\overline{#1}}

\def\2#1{\widetilde{#1}}

\def\3#1{\widehat{#1}}

\def\4#1{\mathbb{#1}}

\def\be{\begin{Exa}}
\def\ee{\end{Exa}}

\def\5#1{\frak{#1}}

\def\6#1{{\mathcal{#1}}}

\def\C{{\4C}}

\def\R{{\4R}}

\def\sumJ{\underset{|J|=k}{{\sum}'}}

\def\sumjq{\underset {j\leq q}\sum}

\def\T{\text}
\newcommand{\Om}{\Omega}
\newcommand{\om}{\omega}

\newcommand{\bom}{\bar{\omega}}

\newcommand{\no}[1]{\|{#1}\|}

\def\NO#1{||#1||^2}
\def\R{{\Bbb R}}

\def\C{{\Bbb C}}

\def\di{\partial}
\def\dib{\bar\partial}
\def\Label#1{\label{#1}}

\def\SPsi{\mathcal S(\Psi)}

\def\simleq{\underset\sim<}

\def\simgeq{\underset\sim>}

\def\T{\text}

\def\1#1{\overline{#1}}

\def\2#1{\widetilde{#1}}

\def\3#1{\widehat{#1}}

\def\4#1{\mathbb{#1}}

\def\5#1{\frak{#1}}

\def\6#1{{\mathcal{#1}}}

\def\C{{\4C}}

\def\R{{\4R}}

\def\sumJ{\underset{|J|=k}{{\sum}'}}

\def\sumjq{\underset{j=1}{\overset{q_0}\sum}}

\numberwithin{equation}{section}

\def\T{\text}

\theoremstyle{plain}

\newtheorem{theorem}{Theorem}[section]

\newtheorem{lemma}[theorem]{Lemma}

\newtheorem{proposition}[theorem]{Proposition}

\theoremstyle{definition}

\newtheorem{definition}[theorem]{Definition}

\theoremstyle{remark}

\newtheorem{remark}[theorem]{Remark}

\begin{document}

\title[Local regularity of the Green operator...]{ Local regularity of the Green operator in a CR manifold of general ``type"}         
\author[L.~Baracco, T.V.~Khanh, S.~Pinton and G.~Zampieri]{Luca Baracco, Tran Vu Khanh, Stefano Pinton and  Giuseppe Zampieri}
\address{Dipartimento di Matematica, Universit\`a di Padova, via 
Trieste 63, 35121 Padova, Italy}
\email{baracco@math.unipd.it,khanhpinton@math.unipd.it,  
zampieri@math.unipd.it}
\maketitle

\begin{abstract}
It is here proved that if a pseudoconvex CR manifold $M$ of hypersurface type has a certain ``type", that we quantify by a vanishing rate $F$ at a submanifold of CR dimension $0$, then $\Box_b$ ``gains $f^2$ derivatives" where $f$ is defined by inversion of $F$. 
Next a general tangential
estimate, ``twisted" by a pseudodifferential operator $\Psi$ is established. The combination of the two yields a general ``$f$-estimate" twisted by $\Psi$, that is, \eqref{1.2} below. 
We apply the twisted estimate  for $\Psi$ which is the composition of a cut-off $\eta$ with a differentiation of order $s$ such as $R^s$ of Section~\ref{s3}. Under the assumption
that $[\di_b,\eta]$ and $[\di_b,[\dib_b,\eta]]$ are superlogarithmic multipliers in a sense inspired to Kohn, we get the local regularity of the Green operator $G=\Box_b^{-1}$.  In particular, if $M$ has ``infraexponential type" along $S\setminus\Gamma$ where $S$ is a manifold of CR dimension $0$ and $\Gamma$ a curve transversal to $T^\C M$, then we have local regularity of $G$. This gives an immediate proof of \cite{BKZ14} in tangential version and of \cite{K00}. The conclusion extends to ``block decomposed" domains for whose blocks the above hypotheses hold separately.
\end{abstract}
\tableofcontents 
\section{Introduction}
\Label{s1}
It has been proved in \cite{KZ10} that if the boundary of a pseudoconvex domain of $\C^n$ has geometric ``type $F$", then there is an ``$f$-estimate" for the $\dib$-Neumann problem for $f=F^*(t^{-1})^{-1}$ where $F^*$ is the inverse function to $F$. The converse is also true (cf. \cite{KZ12}), apart from a loss of accuracy in the estimate which is in most cases negligeable. The succesful approach in establishing the equivalence between the $F$-type and the $f$-estimate consists in triangulating through a potential theoretical condition, namely, the ``$f$-property", that is, the existence of a bounded weight whose Levi-form grows with the rate of $f^2$ at the boundary. This generalizes former work by Kohn \cite{K79}, Catlin \cite{C83}, \cite{C87}, McNeal \cite{MN92} et alii.
 What we prove here is that the $F$ type implies the $f$-estimate for the tangential system $\dib_b$; this is a generalization of  Kohn \cite{K02}. 
 In greater detail, let $M\subset \C^n$ be a pseudoconvex manifold of hypersurface type
and $v$ or $u$ a form in $M$ of a certain degree $h$. We use the microlocal   decomposition into wavelets $u=\sum_{k=1}^{+\infty}\Gamma_ku$ (cf. \cite{K02} proof of Theorem~6.1). 
We consider 
a submanifold  $S\subset M$  of CR dimension $0$, and a real function $F$ satisfying $\frac F{d_S^2}\searrow0$ as the distance $d_S$ to $S$ decreases to $0$. We also use the notation Id for the identity of the complex tangent bundle $T^\C M=TM\cap iTM$. We assume that $M$ has  type $F$ along $S$ in a neighborhood $U$ of  point $z_o\in S$ in the sense that  the  Levi form $(c_{ij})$ of $M$ satisfies
$
(c_{ij})\simgeq \frac{F(d_S)}{d_S^2}\,\T{Id}.
$
Then, there is a bounded family of weights $\{\phi^k\}$  by the aid of which we get the estimate
of the $f$-norm by the Levi form $(c_{ij})$ of $M$ and $(\phi^k_{ij})$ of the $\phi^k$'s.
\bt
\Label{t1.1}
Let $M$ have type $F$ along $S$; then 
\begin{equation}
\Label{1.0,1}
\begin{cases}
\begin{split}
\no{f(\Lambda)v}&\simleq \int_M(c_{ij})(\Lambda^{\frac12}v,\overline{\Lambda^{\frac12}v})\,dV+\sum_{k=1}^{+\infty}\int_M(\phi_{ij}^k)(\Gamma_kv,\overline{\Gamma_kv})\,dV
\\
&+\NO{v}_0,\quad \T{for any $v$ of degree $h\in[1,\dim_{CR}(M)]$},
\end{split}
\\
\begin{split}
\no{f(\Lambda)v}&\simleq \int_M\Big(\T{Trace}(c_{ij})\T{Id}-(c_{ij})\Big)(\Lambda^{\frac12}v,\overline{\Lambda^{\frac12}v})\,dV+\sum_{k=1}^{+\infty}\int_M\Big(\T{Trace}(\phi^k_{ij})\T{Id}-(\phi_{ij})\Big)\times
\\&
\times (\Gamma_kv,\overline{\Gamma_kv})\,dV+\NO{v}_0,\quad \T{for any $v$ of degree $h\in[0,\dim_{CR}(M)-1]$}.
\end{split}
\end{cases}
\end{equation}
\et
The proof is the content of Section~\ref{s2} below.
We denote by $u=u^++u^-+u^0$ the microlocal decomposition of $u$ (cf. \cite{K02} Section~2)
and also use the notation $Q^b $ for the energy $Q^b=\NO{\dib_b v}+\NO{\dib_b^* v}$, and $\mathcal H$ for the space of harmonic forms $\mathcal H=\ker\dib_b\cap\ker\dib_b^*$.
We apply the first of \eqref{1.0,1} for $v=u^+$, resp. the second for $v=u^-$, and plug into a basic estimate. We also use the elliptic estimate for $u^0$ and conclude
\bt
\Label{t1.1,5}
We have
\begin{equation}
\Label{f}
\NO{f(\Lambda)u}\simleq Q^b(u,\bar u)+\NO{u}_0,\quad \T{ for any $u$ of degree }h\in[0,\dim_{CR}(M)].
\end{equation}
\et
As it has been already said, \eqref{f} follows from \eqref{1.0,1} for the common range of degrees $h\in[1,\dim_{CR}(M)-1]$. As for the critical top and bottom degrees, we get the estimate for $u\in\mathcal H^\perp$ from the estimate in nearby degree from closed range of $\dib_b$ and $\dib_b^*$ (\cite{K02} proof of Theorem~7.3 p. 237).

Next, we prove a general basic weighted estimate twisted by a pseudodifferential operator $\Psi$, that is, \eqref{3.1} and \eqref{3.1bis} of Theorem~\ref{t3.1} below. 
We have  to mention that our formula is classical (cf. McNeal \cite{MN06}, \cite{S10}) when $\Psi$ is a function. A recent  application, in which $\Psi$ is a family of cut-off, has been given in  \cite{BPZ14} in the problem of the local regularity of the Green operator $G=\Box_b^{-1}$.
We choose a smooth orthonormal basis of $(1,0)$ forms $\om_1,...,\om_{n-1}$, supplement by a purely imaginary form $\gamma$ and denote the dual basis of vector fields by $\di_{\om_1},...,\di_{\om_{n-1}},T$. We define various constants $c^h_{ij}$'s  as the coefficients of the commutator
$
[\di_{\om_i},\dib_{\om_j}]=c_{ij}^nT+\sum_{j=1}^{n-1}c_{ij}^h\di_{\om_h}-\sum_{j=1}^{n-1}\bar c_{ji}^h\dib_{\om_h}$; sometimes, we also write $c_{ij}$ instead of $c_{ij}^n$. 
We use the notation $\Op$ for an operator of order smaller than $\Psi$ whose support is contained in a conical neighborhood of that of $\Psi$.
Combination of the $f$ estimate with the basic twisted estimate yields 
\bt
\Label{t1.2}
Let $M$ have type $F$ along a CR manifold $S$ of CR dimension $0$ at $z_o$ and $U=U_t$ be suitably small.
 For any form $v=u^+\in C^\infty_c(M\cap U)$ of degree $h\in[1,\dim_{CR}(M)-1]$ we have
\begin{equation}
\Label{1.0,3}
\begin{split}
||f(\Lambda)\Psi v||^2&\le \int(c_{ij})(\Psi T^{\frac12} v,\overline{\Psi T^{\frac12}v})\,dV+\sum_k\int(\phi^k)_{ij}(\Gamma_k\Psi v,\overline{\Gamma_k\Psi v})dV +t\NO{\Psi v}_0
\\
&\simleq Q_{\Psi}^b(v,\overline{v})+\NO{[\di_b,\Psi]\contrazione v}_0+\Big|\sum_h\int (c_{ij}^h)([\di_{\om_h},\Psi](v),\overline{\Psi  v})\,dV\Big|
\\
&\qquad+\Big|\int [\di_b,[\dib_b,\Psi^2]](v,\overline{ v})dV\Big|+Q^b_{\Op}(v,\bar v)+\NO{\Op v}_0+\NO{\Psi v}_0.
\end{split}
\end{equation}
Here $Q^b_{\Psi}=\NO{\Psi\dib_bv}+\NO{\Psi\dib_b^*v}$.

\noindent
{\bf (ii)}\hskip0.2cm
The similar equation holds for $u^-$ in degree $[0,\dim_{CR}(M)-1]$ if we replace $(c_{ij})$, $(\phi^k_{ij})$ and $[\di_b,[\dib_b,\Psi^2]]$ by $-(c_{ij})+\sum_j c_{jj}\T{Id}$, $-(\phi^k_{ij})+\sum_j \phi_{jj}\T{Id}$ and $-[\di_b,[\dib_b,\Psi^2]]+\T{Trace}([\di_b,[\dib_b,\Psi^2]])\T{Id}$ respectively.

\noindent
{\bf (iii)} \hskip0.2cm
Taking summation of the estimate for $v=u^+, v=u^-$ together with the elliptic estimate for $v=u^0$, and using the closed range of $\dib_b$ and $\dib^*_b$ for the critical degrees we get for the full $u\in \mathcal H^\perp$ in degree $h\in [0,\dim_{CR}(M)]$
\begin{equation}
\Label{1.2}
\begin{split}
\NO{f(\Lambda)\Psi u}_0&\simleq Q^b_{\Psi}(u,\bar u)+\NO{[\di_b,\Psi]\contrazione u}_0+\Big|\int_M[\di_b,[\dib_b,\Psi^2]](u^+,\overline{u^+})\,dV\,\Big|
\\
& +\,\Big|\sum_h\int (c_{ij}^h)([\di_{\om_h},\Psi](u),\overline{\Psi  u})\,dV\Big|+\Big|\int_M\Big(-[\di_b,[\dib_b,\Psi^2]](u^-,\overline{u^-})
\\
&\qquad+\T{Trace}([\di_b,[\dib_b,\Psi^2]])\T{Id}\Big)(u^-,\overline{u^-})\,dV\Big|+Q^b_{\Op}(u,\bar u)+\NO{\Op u}_0+\NO{\Psi u}_0.
\end{split}
\end{equation}
\et
The proof is just the superposition of the items (i) and (ii)  of Theorem~\ref{t3.1} below. We have indeed, in Theorem~\ref{t3.1} (i) and (ii)  a more general, weighted version of this estimate. 
We give an application of the general twisted estimate   in which $\Psi$ includes a cut-off $\eta$ and a differentiation of arbitrarily high order $s$ (such as $R^s$ of Section~\ref{s4} below). 
To introduce it, we need the notion of superlogarithmic multipliers which are an obvious variant of the subelliptic multipliers (cf. \cite{K02} Definition~8.1). The crucial point in our discussion is that we consider vector multipliers $g=(g_j)$ and also require a more intense property in which energy is replaced by Levi form, that is, for any $\epsilon$, suitable $c_\epsilon$, and for an uniformly bounded family of weights $\{\phi^k\}$
\begin{equation}
\Label{1.0,4}
\NO{\log(\Lambda)g\contrazione v}\simleq\epsilon\Big(\int_M(c_{ij}(\Lambda^{\frac12}v,\overline{\Lambda^{\frac12}v})\,dV+\sum_{k=1}^{+\infty}\int_M(\phi_{ij}^k)(\Gamma_kv,\overline{\Gamma_kv})\,dV\Big)+c_\epsilon\NO{v}_0.
\end{equation}
We also require that the same estimate holds for $(c_{ij})$ and $(\phi^k_{ij})$ replaced by  $-(c_{ij})+\T{Trace}(c_{ij})\T{Id}$ and $-(\phi^k_{ij})+\T{Trace}(\phi^k_{ij})\,\T{Id}$ respectively. 
With this preliminary we have
\bt
\Label{t1.3}
Assume that there is a system of cut-off $\{\eta\}$ at $z_o$ such that $[\dib_b,\eta]$ and $[\di_b,[\dib_b,\eta]]$ are vector and matrix superlogarithmic multipliers respectively,
and $(c_{ij}^h)$ are subelliptic multipliers. 
 Then $G$ is regular at $z_o$.
\et
The proof is found in Section~\ref{s4}. We combine Theorem~\ref{t1.3} with \ref{t1.1}. This gives back the conclusion of \cite{BKZ14} (in a tangential version) which was in turn  a  generalization of \cite{K00}. It also provides a larger class of hypersurfaces for which $G$ is regular. Let $M$ be the ``block decomposed" hypersurface of $\C^n$  defined by $x_n=\sum_{j=1}^m h^{I^j}(z_{I^j},y_n)$ where $z=(z_{I^1},...,z_{I^m},z_n)$ is a decompostion of  coordinates.
\bt
\Label{t1.4}
Assume that
\begin{equation}
\Label{1.4}
\begin{cases}
 \T{ (a) $h^{I^j}$ has infraexponential type along a totally real $S^{I^j}\setminus\Gamma^{I^j}$ where $S^{I^j}$ is}
\\
\quad\T{  totally real  in $\C^{I^j}\times \C_{z_n}$ and $\Gamma^{I^j}$ is a curve of $\C^{I^j}\times\C_{z_n}$ transversal to $\C^{I^j}\times\{0\}$,}
\\
 \T{(b) $h^j_{z_j}$ are superlogarithmic multipliers},
\\
\T{(c) $c_{ij}^h$ are subelliptic multipliers}.
\end{cases}
\end{equation}
Then,  we have local regularity of $G$ at $z_o=0$. 
\et
In case of a single block $x_n=h^{I^1}$ we regain \cite{BPZ14} and \cite{K00}.
The proof is found in Section~\ref{s4} below.

\noindent
{\it Example} 
Let
$$
(i)\quad
x_n=\sum_{j=1}^{n-1}e^{-\frac1{|z_j|^a}}e^{-\frac1{|x_j|^b}}\qquad\T{ for any $a\geq0$ and for $b<1$}.
$$
Then, \eqref{1.4} (a) is obtained starting from $h^j_{z_j\bar z_j}\simgeq \frac{e^{-\frac1{|x_j|^b}}}{|x_j|^2}$, that is, the condition of type $F_j:=e^{-\frac1{|x_j|^b}}$ along $S_j=\R_{y_j}\times\{0\}$. This yields the estimate of the $f$ norm for  $f(t)=
\log^{\frac1b} (t)$; since $\frac1b>1$, this is superlogarithmic. \eqref{1.4} (b) follows from $|h^j_{z_j}|^2\simleq h^j_{z_j\bar z_j}$ which says that the $h^j_{z_j}$'s are not only superlogarithmic, but indeed $\frac12$-subelliptic,  multipliers. 
Finally, (c) follows from $c_{jj}^h\simleq c_{jj}$ (a consequence of the ``rigidity" of $M$) which shows that these constant are $\frac12$ subelliptic multipliers. 

\eqref{1.4} is the ultimate step of a long sequence of criteria of regularity of $G$, not reduceable in one another, described by the hypersurface models below, in which $a>0$ and $0<b<1$,
\begin{itemize}
\item [(ii)] $x_n=\sum_{j=1}^{n-1}e^{-\frac1{|x_j|^b}}$ Kohn \cite{K02},
\item[(iii)] $x_n=e^{-\sum_{j=1}^{n-1}\frac1{|z_j|^a}}$ Kohn \cite{K00},
\item[(iv)] $x_n=e^{-\frac1{\sum_{j=1}^{n-1}|x_j|^a}}\Big(\sum_{j=1}^{n-1}e^{-\frac1{|x_j|^b}}\Big)$ Baracco-Khanh-Zampieri \cite{BKZ14},
\item[(v)] $x_n=\sum_{j=1}^{n-1}e^{-\frac1{|z_j|^a}}$ Baracco-Pinton-Zampieri \cite{BPZ13},
\item[(vi)] $x_n=\sum_{j=1}^{n-1}e^{-\frac1{|z_j|^a}}x_j^a$ Baracco-Pinton-Zampieri \cite{BPZ14}.
\end{itemize}
Thus, the degeneracy in our model (i) comes as the combination of those of (ii) with (v) (or (vi)).

\section{Estimate of the $f$-norm by the Levi form}
\Label{s2}
Let $M$ be a $C^\infty$ CR-manifold of   $\C^n$ of hypersurface-type, $z_o$ a point of $M$, $U$ an open neighborhood of $z_o$.
Our setting being local, we can find a local CR-diffeomeorphism which reduces $M$ to a hypersurface of $TM+iTM$; therefore, it is not restrictive to assume that $M$ is a hypersurface of $\C^n$ from the beginning. 
We choose a smooth orthonormal basis of $(1,0)$ forms $\om_1,...,\om_{n-1}$, supplement by a purely imaginary form $\gamma$ and denote the dual basis of vector fields by $\di_{\om_1},...,\di_{\om_{n-1}},T$.  We also use the notation $\dib_b$ for the tangential CR-system. For a smooth real function $\phi$, we denote by
$(\phi_{ij})$  the matrix of the Levi form $\di_b\dib_b\phi$. Note that $\phi_{ij}$ differs from $\di_{\om_i}\dib_{\om_j}(\phi)$ because of the presence of the derivatives of the coefficients of the forms $\bar\di_{\om_j}$.
Let $(c_{i\bar j})_{i,j=1,...n-1}$ be the Levi-form $d\gamma|_{T^\C {M}}$ where $T^\C {M}=T{M}\cap iT{M}$.

Let $S\subset M$ be a submanifold of CR-dimension $0$, $d_S$ the Euclidean distance to $S$, and $f:\,\R^+\to\R^+$  a smooth monotonic increasing function such that $f\simleq t^{\frac12}$.
We use the notation $a_k$ for the constant $a_k:=f^{-1}(2^k)$ and $S_{a_k}$ for the strip $S_{a_k}:=\{z\in M:\,\,d_S(z)\le a_k\}$.
\bl
\Label{l2.1}
There is an uniformly bounded family of smooth weights $\{\phi^k\}$ with supp $\phi^k\subset S_{2a_k}$ whose Levi-form satisfies
\begin{equation}
\Label{2.1}
\di_b\dib_b\phi^k\simgeq \begin{cases}
f^2(2^k)&\T{on $S_{a_k}$}
\\
-f^2(2^k)&\T{on $S_{2a_k}\setminus S_{a_k}$},
\\
0&\T{ on $M\setminus S_{2a_k}$}.
\end{cases}
\end{equation}
This also readily implies the same inequalities as \eqref{2.1} with $\di_b\dib_b\phi^k$ replaced by $\Big(\T{Trace}(\di_b\dib_b\phi^k)\,\T{Id}-\di_b\dib_b\phi^k\Big)$.
\el
Note that there is no  assumption about the behavior of $M$ at $S$ in this Lemma.

\bpf
Set
\begin{equation}
\Label{2.2} 
\phi^k=c\chi(\frac{d_S(z)}{a_k})\log(\frac{d_S^2(z))}{a_k^2}+1),
\end{equation}
where $c$ is a constant that will be specified later and $\chi\in C^\infty(0,2)$ is a decreasing cut-off function which satisfies
\begin{equation*}
\begin{cases}
\chi\equiv1&\T{on $[0,1]$},
\\
0\le \chi\le 1 &\T{on $[1,\frac32]$},
\\
\chi\equiv0 &\T{on $[\frac32,2]$}.
\end{cases}
\end{equation*}
Remark that 
\begin{equation*}
\begin{split}
\di_b\dib_b d_S^2&=2\di_b d_S\otimes \dib_b d_S+2d_S\di_b\dib_b d_S
\\
&\ge 2\di_b d_S\otimes \dib_b d_S\\
&\simgeq\T{Id},
\end{split}
\end{equation*}
where the last inequality follows from $\dim_{CR}(M)=0$ (with the agreement that Id denotes the identity of $T^\C M$).

 Now, when $\di_b\dib_b$ hits $\log$, we have
\begin{equation}
\Label{2.3}
\begin{split}
\di_b\dib_b \log (\frac{d_S^2(z)}{a_k^2}+1)
&\simgeq\frac{\di_b d_S\otimes\dib_b d_S+d_S\di_b\dib_b d_S}{a_k^2}
\\
&\simgeq \frac{\T{Id}}{a_k^{2}}=f^2(2^k)\,\T{Id}.
\end{split}
\end{equation}
On the other hand, on $S_{a_k}$, the function  $\chi$ is constant and therefore $\di_b\dib_b\phi^k=\di_b\dib_b \log $. Thus \eqref{2.3} yields the first of 
 \eqref{2.1}. When, instead, $\di_b\dib_b$ hits $\chi$, we have
\begin{equation}
\Label{2.3,5}
\begin{split}
\Big|\di_b\dib_b\chi\Big(\frac{d_S(z)}{a_k}\Big)\Big|&\leq|\ddot\chi|\frac{\di_bd_S\otimes \dib_b d_S}{a_k^2}+|\dot\chi|\frac{\di_b\dib_bd_S}{a_k}
\\
&\underset{\T{ since $\dim_{CR}(S)=0$}}\simleq \frac{\T{Id}}{a_k^2}.
\end{split}
\end{equation}
On the other hand, $\log$ stays bounded on $S_{2a_k}$ and therefore $\di_b\dib_b(\chi)\log\simgeq -a_k^{-2} =-f^{2}(2^k)$. Finally, when $\di_b$ and $\dib_b$ hit $\chi$ and $\log$ separately, we get 
\begin{equation}
\Label{3.5,6}
\begin{split}
\Big|2\Re e \di_b\chi\Big(\frac{d_S}{a_k}\Big)\dib_b\log(\frac{d_S^2}{a_k^2}+1)\Big|&\simleq\Big|2\Re e\dot\chi\frac{\di_b d_S}{a_k}\otimes \frac{2a_k^2d_S\dib_bd_S}{2d_S^2a_k^2}\Big|
\\
&\underset{\T{since $d_S\sim a_k$ on supp$\,\dot\chi$}}\simleq \frac{\di_bd_S\otimes \dib_b d_S}{a_k^2}=f^2(2^k)\,\T{Id}.
\end{split}
\end{equation}
Thus, again, $2\Re e\dib_b\chi\dib_b\log\simgeq -f^2(2^k)\,\T{Id}$.

\epf

As we have seen in the proof of Lemma~\ref{l2.1}, when $\dot\chi$ and $\ddot\chi\neq0$, the Levi form of $\phi^k$ can get negative. However, this annoyance can be well behaved by the aid of the Levi form of $M$.  
Let $F$ be a smooth real function such that $\frac{F(d)}{d^2}\searrow 0$ as $d\searrow0$, denote by $F^*$ the inverse to $F$ and define $f(t):=(F^*(\delta))^{-1}$, for $\delta=t^{-1}$. 
Let $f(\Lambda)$ be the tangential pseudodifferential operator with symbol $f$.
 This is defined by introducing a local  straightening  ${M}\simeq \R^{2n-1}\times\{0\}$  for a defining function $r=0$ of $M$, taking  local coordinates $x\in M$, dual coordinates $\xi$ of $x$ and setting
$$
f(\Lambda)(u)=\int \left(e^{ix\xi}f(\sqrt{1+\xi^2})\int e^{-iy\xi} u(y)dy\right)d\xi.
$$
In particular $\Lambda$ is the standard elliptic pseudodifferential operator with symbol $\sqrt{1+\xi^2}$. 
\bd
\label{type}
We say that $M$ has type $F$ along $S$ in  a neighborhood $U$ of $z_o$, if
\begin{equation}
\Label{2.4}
(c_{ij})\simgeq \frac{F(d_S)}{d_S^2}\T{Id}\quad\T{on $U$}.
\end{equation}
\ed
Note that \eqref{2.4} implies
\begin{equation}
\Label{2.4bis}
\Big(\T{Trace}(c_{ij})\T{Id}-(c_{ij})\Big)\simgeq \frac{F(d_S)}{d_S^2}\T{Id}\quad\T{on $U$}.
\end{equation}
\bp
\Label{p2.1}
Let $M$ have type $F$ along S of CR dimension 0. Then 
\begin{equation}
\Label{2.5}
\begin{cases}
\no{f(\Lambda)\Gamma_kv}^2_0\simleq \int_M(c_{ij})(\Gamma_k\Lambda^{\frac12}v,\overline{\Gamma_k\Lambda^{\frac12}v})\,dV+ \int_M (\phi^k_{ij})(\Gamma_kv,\overline{\Gamma_kv})\, dV+\no{\Gamma_kv}_0^2,\,\,h\in[1,n-1],
\\
\begin{split}
\no{f(\Lambda)\Gamma_kv}^2_0&\simleq \int_M\Big(\T{Trace}(c_{ij})\T{Id}-(c_{ij})\Big)(\Gamma_k\Lambda^{\frac12}v,\overline{\Gamma_k\Lambda^{\frac12}v})\,dV
\\
&+ \int_M \Big(\T{Trace}(\phi^k_{ij})\T{Id}-(\phi^k_{ij})\Big)(\Gamma_kv,\overline{\Gamma_kv})\, dV+\no{\Gamma_kv}_0^2,\,\,h\in[0,n-2].
\end{split}
\end{cases}
\end{equation}
\ep
\bpf
We set $a_k=f^{-1}(2^k)=F^*(2^{-k})$, $S_{a_k}=\{z:\,d_S(z)<a_k\}$ and denote by $\lambda(z)$ the minimum of the $n-1$ eigenvalues of $(c_{ij})$ at $z$. 
We start from the first of \eqref{2.5}. We have 
\begin{equation}
\Label{2.7}
\begin{split}
\no{\Gamma_kv}^{M\setminus S_{a_k}}_0&\simleq \underset{z\in M\setminus S_{a_k}}\max \frac{2^{-\frac k2}}{\lambda(z)^{\frac12}}\Big(\no{
\lambda^{\frac12}\Gamma_k\Lambda^{\frac12} v}^{M\setminus S_{a_k}}+\no{\Gamma_kv}_{-\frac12}\Big)
\\
&\simleq \frac{a_k2^{-\frac k2}}{F(a_k)^{\frac12}}\Big(\sqrt{\int_{M\setminus S_{a_k}}(c_{ij})(\Gamma_k\Lambda^{\frac12}v,\overline{\Gamma_k\Lambda^{\frac12}v})\,dV}+2^{-\frac k2}\no{\Gamma_kv}^{M\setminus S_{a_k}}_0\Big)
\\
&\simleq f^{-1}(2^{k})\Big(\sqrt{\int_{M\setminus S_{a_k}}(c_{ij})(\Gamma_k\Lambda^{\frac12}v,\overline{\Gamma_k\Lambda^{\frac12}v})\,dV}+\no{\Gamma_kv}_0\Big).
\end{split}
\end{equation}
Recalling that $f(\Lambda_\xi)\equiv f(2^k)$ on supp$\,\Gamma_k$, this gives 
\begin{equation}
\Label{*}
\no{f(\Lambda)\Gamma_kv}^{M\setminus S_{a_k}}_0\simleq \sqrt{\int_{M\setminus S_{a_k}}(c_{ij})(\Gamma_k\Lambda^{\frac12}v,\overline{\Gamma_k\Lambda^{\frac12}v})\,dV}+\no{\Gamma_kv}_0.
\end{equation}
Now, on $S_{2a_k}\setminus S_{a_k}$, $(\phi^k)_{ij}$ can get negative. However, using the second of \eqref{2.1} and tuning the choice of $c$, independent of $k$,  in the definition of $\phi^k$ so that $2^k(c_{ij})+(\phi^k)_{ij}\geq \frac{f^2(2^k)}2$ on $S_{2a_k}\setminus S_{a_k}$, we have that not only \eqref{*} but also \eqref{2.5} holds on $M\setminus S_{a_k}$. 

Finally, on $S_{a_k}$, $(\phi^k)_{ij}$ satisfies the first of \eqref{2.1} and therefore
$$
\no{f(\Lambda)\Gamma_kv}^{ S_{a_k}}_0\simleq \sqrt{\int_{ S_{a_k}}(\phi^k_{ij})(\Gamma_kv,\overline{\Gamma_kv})\,dV}+\no{\Gamma_kv}_0.
$$
This shows how \eqref{2.5} follows from \eqref{2.4}. In the same way we can see that the second follows from \eqref{2.4bis}.
\epf

\noindent
{\it Proof of Theorems~\ref{t1.1} and \ref{t1.1,5}.} The proof  of \eqref{1.0,1} just consists in taking summation over $k$ in \eqref{2.5}. As for \eqref{f} in degrees $h\in [1,n-2]$, it follows from the combination of the first (resp. the second) of \eqref{1.0,1} for $v=u^+$ (resp. $v=u^-$), in addition to the classical basic tangential estimates and  the elliptic estimate for $u^0$. As for the critical degree $h=0$ and $h=n-1$ in \eqref{f}, it follows from writing $u=\dib_b^*w$ and $u=\dib_bw$ respectively (by closed range) and by using the estimate already established for $w$ in the non-critical degrees $1$ and $n-2$ respectively.

\section{The tangential H\"ormander-Kohn-Morrey formula twisted by a pseudodifferential operator}
\Label{s3}
Let $M$ be a CR manifold of hypersurface type of $\C^n$, $\dib_b$ the tangential Cauchy-Riemann system, $\dib_b^*$ the adjoint system. Our discussion is local and we can therefore assume that $M$ is in fact a hypersurface. For a neighborhood $U$ of a point $z_o\in M$, we identify $U\cap M$ to $\R^{2n-1}$ with coordinates $x$ and dual coordinates $\xi$, and consider a pseudodifferential operator $\Psi$ with symbol $\mathcal S(\Psi)(x,\xi)$. For notational convenience we assume that the symbol is real. We also use the notation $L^2_\phi$ for the $L^2$ space weighted by $e^{-\phi}$, $Q^b=\no{\dib_b u}^2+\no{\dib^*_b u}^2$ for the energy, and $Q^{b\,\phi}_\Psi=\no{\Psi\dib_b u}^2_\phi+\no{\Psi\dib^*_b u}^2_\phi$ for the  energy weighted by $\phi$ and twisted by $\Psi$. We consider the pseudodifferential decomposition of the identity by Kohn $\T{Id}=\Phi^++\Phi^-+\Phi^0$ modulo $\T{Op}^{-\infty}$. We consider a basis of $(1,0)$ forms $\om_1,...,\om_{n-1}$ the conjugate basis $\bom_1,...,\bom_{n-1}$ and complete by a purely imaginary form $\gamma$. We denote by  $\di_{\om_1},...,\di_{\om_{n-1}},\dib_{\om_1},...,\dib_{\om_{n-1}},T$ 
the dual basis of vector fields.
$M$ being a hypersurface defined, say, by $r=0$, we can supplement the $\om_j$'s to a full basis of $(1,0)$ forms in $\C^n$ by adding $\om_n=\di r$. Then $\gamma=\om_n-\bar\om_n$ and $T=\di_{\om_n}-\di_{\bar\om_n}$. We describe the commutators by
\begin{equation}
\Label{commutator}
\begin{split}
[\di_{\om_i},\di_{\bar\om_j}]&=\sum_{j=1}^{n}c_{ij}^h\di_{\om_h}-\sum_{j=1}^{n}\bar c_{ji}^h\dib_{\om_h}
\\
&=c^n_{ij}T+,\sum_{j=1}^{n-1}c_{ij}^h\di_{\om_h}-\sum_{j=1}^{n-1}\bar c_{ji}^h\dib_{\om_h};
\end{split}
\end{equation}
We also write $c_{ij}$ instead of $c^n_{ij}$. 

 For a cut-off $\eta\in C^\infty_c(U\cap M)$ we write $u^+:=\eta\Phi^+u,\,\,u^-=\eta\Phi^-u,\,\,u^0=\eta\Phi^0u,\,\,T^{\overset\pm0}=\eta T\Phi^{\overset\pm0}$. We note that $\mathcal S(T)>0$ on supp$\,\mathcal S(\Phi^+)$ (resp. $\mathcal S(T^-)>0$ on supp$\,\mathcal S(\Phi^-)$) and therefore $T^{\frac12}$ (resp. $(T^-)^{\frac12}$) makes sense when acting on $u^+$ (resp. $u^-$).  
We make the relevant remark that
\begin{equation*}
\begin{cases}
\mathcal S(T)\sim \Lambda\T{ on supp$\,\mathcal S(\Phi^+)$}, \qquad \mathcal S(T^-)\sim \Lambda\T{ on supp$\,\mathcal S(\Phi^-)$},
\\
\{\mathcal S(\di_{\om_j}\}_{j=1,...,n-1}\sim \Lambda\T{ and } \mathcal S(\dib_{\om_j}\}_{j=1,...,n-1}\sim \Lambda\quad\T{on supp$\,\mathcal S(\Phi^0)$}.
\end{cases}
\end{equation*}
We denote by $\Op$ , resp. $\T{Op}^0$, an operator of order $2\T{ord}(\Psi)-\frac12$, resp. $0$, whose support is contained in supp$\,\Psi$; we also assume that $\T{Op}^0$ only depends on the $C^2$-norm of $M$ and, in particular, is independent of $\phi$ and $\Psi$.
\bt
\Label{t3.1}
(i) We have for every smooth form $v=u^+$ of degree $h\in[1,n-1]$ 
\begin{equation}
\Label{3.1}
\begin{split}
\int_M&e^{-\phi}(c_{ij})(T^{\frac12}\Psi v,\overline{T^{\frac12}\Psi v})dV+\int_Me^{-\phi}\Big((\phi_{ij})-\frac12(c_{ij})T(\phi)\Big)(\Psi v,\overline{\Psi  v})dV+\NO{\Psi \bar\nabla v}_\phi
\\
&\simleq Q^{b\,\phi}_{\Psi}(v,\overline{v})+\NO{[\di_b,\Psi]\contrazione v}_\phi+\NO{[\di_b,\phi]\contrazione \Psi v}_\phi+\Big|\sum_{h=1}^{n-1}\int (c_{ij}^h)([\di_{\om_h},\Psi](v),\overline{\Psi  v})\,dV\Big|
\\
&+\Big|\int_M e^{-\phi}[\di_b,[\dib_b,\Psi^2]](v,\overline{v})dV\Big|+Q^{b\,\phi}_{\Op}(v,\bar v)+\NO{\Op v}_\phi+\no{\Psi v}^2_\phi.
\end{split}
\end{equation}
 Here we are using the notation $Q^{b\,\phi}_{\Psi}=\NO{\Psi\dib_b v}_\phi+\NO{\Psi\dib_b^* v}_\phi$.

\noindent (ii) We also have, for $v=u^-$ smooth of degree $h\in [0,n-2]$
\begin{equation}
\Label{3.1bis}
\begin{split}
\int_M&e^{-\phi}\Big(-(c_{ij})+\sum_jc_{jj}\T{Id}\Big)((T^-)^{\frac12}\Psi v,\overline{(T^-)^{\frac12}\Psi v})+\NO{\Psi \nabla v}_\phi
\\
&+\int_Me^{-\phi}\left(\left(-(\phi_{ij})+\sum_j\phi_{jj}\T{Id}\right)+\frac12\left((c_{ij})T(\phi))-(\sum_jc_{jj})T(\phi)\right)\right)(\Psi v,\overline{\Psi v})dV
\\
&\simleq Q^{b\,\phi}_\Psi(v,\overline{v})+\NO{[\di_b,\Psi]\contrazione v}_\phi+\NO{[\di_b,\phi]\contrazione \Psi v}_\phi+\Big|\sum_{h=1}^{n-1}\int\Big(-(c_{ij}^h)+\sum_jc_{jj}^h\T{Id}\Big)([\di_{\om_h},\Psi](v),\times
\\
&\times\overline{\Psi  v})\,dV\Big|+\Big|\int_M e^{-\phi}\Big(-[\di_b,[\dib_b,\Psi^2]]+\T{Trace}([\di_b,[\dib_b,\Psi^2]])\T{Id}\Big)(v,\overline{v})dV\Big|
\\
&+Q^{b\,\phi}_{\Op}(v,\bar v)+\NO{\Op v}_\phi+\no{\Psi v}^2_\phi.
\end{split}
\end{equation}
\et
Clearly $u^0$ is subject to elliptic estimates. These, combined with \eqref{3.1}, \eqref{3.1bis} yield an estimate for the full $u$ in degrees $[1,n-2]$ and then also for $u\in\mathcal H^\perp$ in degree $k\in[0,n-1]$ by closed range.

\br
The formula also holds for $\Psi$ complex: in this case one replaces $\Psi^2$ by $|\Psi|^2$ and add the additional error term $[\di_b,\bar\Psi]\contrazione$ to the already existing $[\di_b,\Psi]\contrazione$.
\er

\bpf
We start from
\begin{equation}
\Label{35}
\begin{split}
\di_b\dib_b\phi&=\di_b(\sum_j\bar \di_{\om_j}(\phi)\bom_j)
\\&\sum_{ij}\left(\di_{\om_i}\dib_{\om_j}(\phi)+\sum_h\overline{c_{ji}^h} \dib_{\om_h}(\phi)\right)\om_i\wedge\bom_j.
\end{split}
\end{equation}
Similarly,
\begin{equation}
\Label{36}
\begin{split}
\dib_b\di_b\phi&=\dib_b(\sum_j \di_{\om_j}(\phi)\om_j)
\\&=\sum_{ij}\left(-\dib_{\om_j}\di_{\om_i}(\phi)-\sum_h{c_{ij}^h}\di_{\om_h}(\phi)\right)\om_i\wedge\bom_j.
\end{split}
\end{equation}
Differently from the ambient $\dib$-system on $\C^n$, we do not have  $\di_b\dib_b=\dib_b\di_b$ and in fact, 
combining \eqref{35} with \eqref{36}, we can describe $(\phi^b_{ij})$, the matrix of $\frac12(\di_b\dib_b-\dib_b\di_b)(\phi)$, by 
\begin{equation}
\Label{3.3}
\begin{split}
\phi_{ij}^b&=\langle\frac12(\di_b\dib_b-\dib_b\di_b)(\phi),\di_{\om_i}\wedge\dib_{\om_j}\rangle
\\
&\underset{\T{by \eqref{35}, \eqref{36}}}=\frac12\Big(\Big(\di_{\om_i}\dib_{\om_j}+\dib_{\om_j}\di_{\om_i}\Big)(\phi)+\sum_{h=1}^{n-1}\bar c_{ji}^h\dib_{\om_h}(\phi)+c_{ij}^h\di_{\om_h}(\phi)\Big)
\\
&=\dib_{\om_j}\di_{\om_i}(\phi)+\frac12\Big([\di_{\om_i},\dib_{\om_j}](\phi)+\sum_h\bar c_{ji}^h\dib_{\om_h}(\phi)+\sum_hc_{ij}^h\di_{\om_h}(\phi)\Big)
\\
&\underset{\T{\eqref{commutator}}}=\dib_{\om_j}\di_{\om_i}(\phi)+\frac12c_{ij}T(\phi)+\sum_hc^h_{ij}\di_{\om_h}(\phi).
\end{split}
\end{equation}
We consider now
\begin{equation}
\Label{3.3bis}
e^{\phi}\Psi^{-2}[\dib_{\om_i},e^{-\phi}\Psi^2]=-\phi_{\bar \om_i}+2\frac{[\dib_{\om_i},\Psi]}\Psi+\frac{\Opbis}{\Psi^2},
\end{equation}
whose sense is fully clear when both sides are multiplied by $\Psi^2$.
In other terms, we have
\begin{equation}
\Label{3.3quater}
\dib^*_{e^{-\phi}\Psi^2}=\dib^*+\di\phi\contrazione-2\frac{[\di,\Psi]}\Psi\contrazione+\frac{\Opbis}{\Psi^2}+\T{Op}^0.
\end{equation}
This leads us to define the transposed operator $\delta_{\om_i}$ to $\dib_{\om_i}$ by
\begin{equation}
\Label{3.3ter}
\delta_{\om_i}:=\di_{\om_i}-\phi_{\om_i}+2\frac{[\di_{\om_i},\Psi]}\Psi+\frac{\Opbis}{\Psi^2}+\T{Op}^0.
\end{equation}
With these preliminaries we have
\begin{equation}
\Label{3.4}
\begin{split}
[\delta_{\om_i},\bar\di_{\om_j}]&=c_{ij}T+\sumh c_{ij}^h\delta_{\om_h}-\sumh\bar c_{ji}^h\dib_{\om_h}+\Big(\phi_{ij}^b-\frac12c_{ij}T(\phi)\Big)
\\&-2\sum_h c_{ij}^h\frac{[\di_{\om_h},\Psi]}{\Psi}+\frac{[\di_{\om_i},[\dib_{\om_j},\Psi]]}\Psi+\frac{[\di_{\om_i},\Psi]\otimes[\dib_{\om_j},\Psi]}{\Psi^2}+\frac{\Opbis}{\Psi^2}+\T{Op}^0.
\end{split}
\end{equation}
We remember now that there are two equally reasonable definition of the pseudodifferential action
\begin{equation}
\Label{pseudo}
\Psi(w)=\begin{cases}
(i)\hskip0.2cm\int e^{ix\xi}\SPsi(x,\xi)\tilde w(\xi)d\xi
\\
(ii)\hskip0.2cm \int e^{ix\xi}(\widetilde{\SPsi}(\cdot,\xi)*\tilde w)d\xi,
\end{cases}
\end{equation}
where $\tilde w$ denotes the Fourier transform.
Up to  error terms of type $\Op$, we have
\begin{equation*}
\begin{split}
\NO{\Psi(w)}&\sim(\Psi w, \Psi w)
\\
&\underset{\T{Plancherel and \eqref{pseudo} (ii)}}\sim(\widetilde{\Psi(w)},\widetilde{\SPsi}(\cdot,\xi)*\tilde w)
\\
&\sim\int\widetilde {\Psi(w)}(\xi)\overline{\int\widetilde{\SPsi}(\xi-\eta,\xi)\tilde w(\eta)d\eta}d\xi
\\
&\underset{\T{$\overline{\widetilde{\SPsi}}(\xi-\eta,\xi)\sim\widetilde{\overline{\SPsi}}(\eta-\xi,\xi)$}}=\int\Big(\int\widetilde{\Psi(w)}(\xi)\widetilde{\overline{\SPsi}}(\eta-\xi,\xi)d\xi\Big)\bar{\tilde w}(\eta)d\eta
\\
&\underset{\T{\eqref{pseudo} (i)}}\sim\int\widetilde{\bar\Psi\Psi(w)}(\eta)\bar{\tilde w}(\eta)d\eta
\\
&\underset{\T{Plancherel}}\sim(|\Psi|^2w,w).
\end{split}
\end{equation*}
For the same reason $(\Psi^2w,w)\sim\int |\Psi|^2|w|^2\,dV$ and therefore
$$
\NO{\Psi(w)}\sim \int|\Psi|^2|w|^2\,dV.
$$
Adding the weight $\phi$ and recalling that in our discussion $\Psi$ is  real,
\begin{equation}
\Label{reasonable}
\NO{\Psi\dib_b^{(*)}v}_\phi=\int e^{-\phi }\Psi^2|\dib_b^{(*)}v|^2\,dV+\NO{\Op(\dib^{(*)}v)}_\phi,
\end{equation}
where $\dib_b^{(*)}$ denotes either $\dib_b$ or $\dib_b^*$.
We are ready for the proof of  \eqref{3.1}; we prove it only for $v=u^+$, the proof of \eqref{3.1bis} for $v=u^-$ being similar. We have
\begin{multline}
\int_\Om e^{-\phi}(c_{ij})(T\Psi v,\overline{\Psi v})+\int_\Om[\di_b,[\dib_b,e^{-\phi}\Psi^2]](v,\overline{v})dV
\\
-\NO{[\di_b,\phi]\contrazione \Psi v}_\phi-\NO{[\di_b,\Psi]\contrazione  v}_\phi
+\NO{\Psi\bar\nabla v}_\phi
\\
\simleq \NO{\Psi\dib_b v}_\phi+\NO{\Psi(\dib_b)^*_{e^{-\phi}\Psi^2}v}_\phi+sc\NO{\Psi\bar\nabla v}_\phi+\Big|\sum_h\int_\Om e^{-\phi}(c_{ij}^h)([\di_{\om_h},\Psi]v,\overline{\Psi v})\,dV\Big|
\\
+Q^{b\,\phi}_{\Op}(v,v)+\NO{\Op v}_\phi+\NO{\Psi v}_\phi,
\end{multline}
or, according to \eqref{3.4} and after absorbing the term which comes with sc,
\begin{multline}
\Label{3.6}
\int_{M}e^{-\phi} c_{ij}(T\Psi v,\overline{\Psi v})dV+\int_{M}e^{-\phi} \phi_{ij}(\Psi v,\overline{\Psi v})dV-\NO{[\di_b,\phi]\contrazione \Psi v}_\phi
\\
+\int_M e^{-\phi}[\di_i,[\dib_j,\Psi^2]](v,\overline{v})dV
-\NO{[\di_b,\Psi]\contrazione v}_\phi+\NO{\Psi\bar\nabla  v}_\phi
\\
\simleq \NO{\Psi\dib_b v}_\phi+\NO{\Psi(\dib_b)^*_{e^{-\phi}\Psi^2}v}_\phi
+\Big|\sum_h\int_\Om e^{-\phi}(c_{ij}^h)([\di_{\om_h},\Psi]v,\overline{\Psi v})\,dV\Big|
\\+Q_{\Op}(v,v)
+\NO{\Op v}_\phi+\NO{\Psi v}_\phi.
\end{multline}
To carry out our proof we need to replace $(\dib_b)^*_{e^{-\phi}\Psi^2}$ by $\dib_b^*$. We have from \eqref{3.4}
\begin{multline}
\Label{3.7}
\NO{\Psi (\dib_b)^*_{e^{-\phi}\Psi^2}v}_\phi\leq \NO{\Psi\dib_b^* v}_\phi+\NO{\Psi\di \phi\contrazione \Psi^2v}_\phi+\NO{[\di_b,\Psi]\contrazione v}_\phi+\NO{\Op v}_\phi
\\
+\underset{\T{\#}}{\underbrace{2\Big|\Re e (\Psi\dib_b^* v,\overline{\Psi\di_b\phi\contrazione v})_\phi\Big|+2\Big|\Re e (\Psi\dib_b^* v,\overline{[\di_b,\Psi]\contrazione v})_\phi+2\Big|\Re e (\Psi\di_b\phi\contrazione v,\overline{[\di_b,\Psi]\contrazione v})_\phi\Big|}}.
\end{multline}
We next estimate by Cauchy-Schwarz inequality
$$
\#\simleq \NO{\Psi\dib_b^* v}_\phi+\NO{\Psi\di_b\phi\contrazione v}_\phi+\NO{[\di_b,\Psi]\contrazione v}_\phi.
$$
We move the third, forth and fifth terms from the left to the right of \eqref{3.6}, and get \eqref{3.1} with $(T\Psi v,\Psi v)$ instead of $(T^{\frac12}\Psi v,T^{\frac12}\Psi v)$. But they only differ for
$$
\Big|\int_M e^{-\phi}\Big ([(c_{ij},T^{\frac12}](T^{\frac12}\Psi v,\Psi v)\Big)\,dV\Big|\simleq \NO{\Psi v}_0,
$$
which is negligeable. 

\epf

We go back to the family of weights of Theorem~\ref{t1.1} and Proposition~\ref{p2.1}. We apply \eqref{3.1} (resp. \eqref{3.1bis}) for $\phi=\phi^k+t|z'|^2$ (resp. $\phi=\phi^k-t|z'|^2$). First, we note that they are absolutely uniformly bounded with respect to $k$; they can be made bounded in $t$ by taking $U=\{z:\,|z'|<\frac1t\}$. (In particular, by boundedness, they can be removed from the norms.) Possibly by raising to exponential, boundednes implies ``selfboundedness of the gradient" when $\phi$ is plurisubharmonic. In our case, in which to be positive is not $(\phi^k_{ij})$ itself but $2^k(c_{ij})+(\phi^k_{ij})$, we have, for $|z'|$ small
\begin{equation}
\Label{supernova}
\begin{split}
|\di_b\phi|^2&=|\di_b(\phi^k+t|z'|^2)|^2
\\
&\simleq |\di_b\phi^k|^2+t^2|z|^2
\\
&\leq 2^k(c_{ij})+(\phi^k_{ij})+t.
\end{split}
\end{equation}
 So $\no{\di_b\phi\contrazione  \Psi u^\pm}^2$ can be removed from the right side of both \eqref{3.1} and \eqref{3.1bis}. 
Also, the term $-\frac12(c_{ij})T(\phi)(v,\bar v)$ is controlled by $(c_{ij})(T^{\frac12}v,\overline{T^{\frac12}v})$ by Sobolev interpolation. 
We then combine Proposition~\ref{p2.1} with Theorem~\ref{t3.1} formula \eqref{3.1} for the weight $\phi^k+t|z'|^2$ (resp. formula \eqref{3.4} for the weight $(\phi^k-t|z'|^2$), and notice that $T^{\frac12}\sim\Lambda^{\frac12}$ on supp$\,\Psi^+$ (resp. $(T^-)^{\frac12}\sim\Lambda^{\frac12}$ on supp$\,\Psi^-$). 
Also, on the right of \eqref{3.1} and \eqref{3.1bis},  one reduces $\NO{\Op v}_\phi$ to $\NO{v}_\phi$ by induction and estimates all terms $Q^\phi_{\Op}$ and $Q^\phi_{{\T{Op}}^{\T{ord}(\Psi)-j}}\,\,j\ge1$ by a common $Q^\phi_{\Psi'}$.

\noindent
{\it Proof of Theorem~\ref{t1.2}.}\hskip0.2cm We have to use \eqref{3.1} with the above choice of the weight $\phi$ and take summation over $k$; this yields \eqref{1.0,3} for $v=u^+$. The twin estimate for $v=u^-$ follows from \eqref{3.1bis} by similar procedure. Finally, \eqref{1.2} comes as the combination of \eqref{1.0,3} for $v=u^+$, the twin for $v=u^-$ and the elliptic estimate for $v=u^0$.

\hskip13cm$\Box$

\section{A criterion of hypoellipticity of the Kohn Laplacian}
\Label{s4}
Let  $M$ be a pseudoconvex, hypersurface type manifold of $\C^n$, $\Box_b=\dib_b\dib^*_b+\dib^*_b\dib_b$   the Kohn Laplacian of $M$, and $G:=\Box_b^{-1}$ the Green operator.  

{\it \bf Proof of Theorem~\ref{t1.3}}\hskip0.2cm
Our program is to prove that for any cut-off $\eta_o\in C^\infty_c(U)$ with $\eta_o\equiv1$ in a neighborhood of $z_o$, for suitable $\eta\succ\eta_o$, that is $\eta|_{\T{supp}\,\eta_o}\equiv1$, for any $s$ and suitable $U$, we have
\begin{multline}
\Label{4.1}
\no{\eta_o u}_s\simleq \no{\eta\dib_b u}_s+\no{\eta\dib^*_b u}_s+\no{u}_0\quad \T{for any $u\in \mathcal H^\perp\cap C^\infty(M\cap U)$}
\\
\T{  in any degree $k\in[0,n-1]$}.
\end{multline}
If we are able to prove \eqref{4.1}, we have immediately the exact local $H^s$-regularity of $\dib^*_bG$ and $\dib G$ over $\ker\dib$ and $\ker\dib^*$ respectively. 
From this, we get the (non-exact) regularity of the Szeg\"o $S=\T{Id}-\dib_b^*G\dib_b$ and anti-Szeg\"o $S^*=\T{Id}-\dib_bG\dib_b^*$ projection respectively. 
(At this stage we need to apply the method of the elliptic regularization to pass from $C^\infty$- to $H^s$-forms.)
From this the (non-exact) regularity of $G$ itself follows (cf. e.g. the proof of Theorem~2.1 of \cite{BKZ14}).
 Along with $\eta_o\prec\eta$, we consider an additional 
 cut-off $\sigma$ with $\eta_o\prec\sigma\prec\eta$ and denote by $R^s$ the pseudodifferential  operator with symbol $(1+|\xi|^2)^{\frac {s\sigma(a)}2}$.  
According to Proposition~2.1 of \cite{BKZ14}, there is no restriction on the degree of $u$; thus $u$ can be either a form or a function.
By Section~\ref{s3} above, we can prove \eqref{4.1} separately on each term of the microlocal decomposition of $u=u^++u^-+u^0$; since $u^0$ has elliptic estimate and $u^-$ can be reduced to $u^+$ by star-Hodge correspondence, we prove the result only  for $v=u^+$. We start from
\begin{equation}
\Label{4.11}
\begin{split}
\no{\Lambda^s\eta_o v}&\simleq \no{R^s\eta_o v}+\no{v}
\\
&=\no{R^s\eta_o\eta^2v}+\no{v}
\\
&\le \no{R^s\eta^2v}+\no{[R^s,\eta_o]\eta^2v}+\no{v}
\\
&\simleq\no{R^s\eta^2v}+\no{v}
\\
&\simleq \no{\eta R^s\eta v}+\no{[R^s,\eta]\eta v}+\no{v}
\\
&\simleq \no{\eta R^s\eta v}+\no{v},
\end{split}
\end{equation} 
(cf. \cite{K02} Section~7).
Next, we apply Theorem~\ref{t3.1}  for $\Psi=\eta R^s\eta$, 
What we have to describe are the error terms in the right of \eqref{3.1}, \eqref{3.1bis}, that is, $[\di_b,\eta R^s\eta]$ and $[\di_b,[\dib_b,\eta R^s\eta]]$. Since the argument is similar for the two, we only treat the first. We have by Jacobi identity
\begin{equation}
\Label{4.6}
\begin{split}
[\di_b,\eta R^s\eta]&=[\di_b,\eta]R^s\eta+\eta[\di_b,R^s]\eta+\eta R^s[\di_b, \eta]
\\
&=[\di_b,R^s]+\T{Op}^{-\infty}.
\end{split}
\end{equation}
In fact, since supp$\,\di_b\eta\cap\,\T{supp}\,\sigma=\emptyset$, then the first and last terms in the right of the first line of \eqref{4.6} are operators of order $-\infty$ and can therefore be disregarded. As for the central term, we have
\begin{equation}
\Label{**}
[\di_b,R^s]=\di_b(\sigma)\log(\Lambda)R^s.
\end{equation}
Now, our hypothesis is that
\begin{multline}
\Label{4.7}
\no{\log(\Lambda)\di_b\sigma\contrazione \eta R^s\eta v}^2\le\epsilon\Big(\int_M(c_{ij})(\Lambda^{\frac12}\eta R^s\eta v,\overline{\Lambda^{\frac12}\eta R^s\eta v})\,dV
\\
+\sum_{k=1}^{+\infty}\int_M \Big((\phi^k_{ij})(\eta R^s\eta \Gamma_kv,\overline{\eta R^s\eta \Gamma_kv})\Big)\,dV\Big)
+c_\epsilon\NO{\eta R^s\eta v}.
\end{multline}
Altogether, we get
\begin{equation}
\Label{4.8}
\begin{split}
t\NO{\Lambda^s\eta_ov}&\underset{\T{\eqref{4.11}}}\simleq t\NO{\eta R^s\eta v}_0+\NO{v}_0
\\
&\simleq \Big(\int_M(c_{ij})(\Lambda^{\frac12}\Psi v,\overline{\Lambda^{\frac12}\Psi v})\,dV+\sum_{k=1}^{+\infty}\int_M(\phi_{ij}^k)(\Gamma_k\Psi v,\overline{\Gamma_k\Psi v})\,dV\Big)+t\NO{\eta R^s\eta v}_0+c_\epsilon\NO{v}_0
\\
&\underset{\T{by the second of \eqref{1.0,3}}}\simleq Q^b_{\eta R^s\eta}(v,\bar v)+\NO{[\di_b,\eta R^s\eta]\contrazione v}_0+\Big|\int_M[\di_b,[\dib_b,\eta R^s\eta]](v,\bar v)\,dV\Big|
\\&\hskip0.3cm+\Big|\sum_h\int (c_{ij}^h)([\di_{\om_h},\Psi](v),\overline{\Psi  v})\,dV\Big|+Q_{\Op}^b(v,\bar v)+\NO{\Op v}_0
\\
&
\underset{\T{\eqref{**} and \eqref{1.4} (c)}}\simleq  Q^b_{\eta R^s\eta}(v,\bar v)+\NO{\di_b(\sigma)\log(\Lambda)\eta R^s\eta v}_0+Q_{\Op}^b(v,\bar v)+\NO{\eta'v}_{s-\epsilon}
\\
&\underset{\T{\eqref{4.7}}}\simleq Q^b_{\eta R^s\eta}(v,\bar v)+\epsilon\Big(\int_M(c_{ij})(\Lambda^{\frac12}\eta R^s\eta v,\overline{\Lambda^{\frac12}\eta R^s\eta v})\,dV+\sum_k\int \Big((\phi^k_{ij})(\eta R^s\eta \Gamma_kv,\times
\\
&\times \overline{\eta R^s\eta \Gamma_kv})\Big)\,dV\Big)
+c_\epsilon\NO{\eta R^s\eta v}_0+Q_{\Op}^b(v,\bar v)+\NO{\eta'v}_{s-\epsilon}
\\
&\underset{\T{absorbtion in the second line}}\simleq Q^b_{\eta R^s\eta}(v,\bar v)+Q_{\Op}^b(v,\bar v)+\NO{\eta R^s\eta v}_0+\NO{\eta'v}_{s-\epsilon}
\\
&\underset{\T{absortion by means of $t$}}\simleq Q^b_{\eta R^s\eta}(v,\bar v)+Q_{\Op}^b(v,\bar v)+\NO{\eta'v}_{s-\epsilon}.
\end{split}
\end{equation}
Now, the $s-\epsilon$ norm is reduced to $0$ norm by induction over $j$ with $j\epsilon>s$, and $Q_\eta R^s\eta$ and the various $Q_{\T{Op}^{s-\epsilon j-1}}$ are estimated by a common  $Q_{\eta'\Lambda^s}$. In conclusion, we have got \eqref{4.1} with the notational difference of $\eta'$ instead of $\eta$. 

\hskip13cm $\Box$

\noindent
{\it Proof of Theorem~\ref{t1.4}.} \hskip0.2cm 
We choose our cut-off starting from a cut-off $\chi$ in one real variable and setting $\eta=\Pi_j\chi(|z_{I^j}|)\chi(|y_n|)$. We have
\begin{itemize}
\item[(a)] supp$\,\di_{z_{I^j}}\chi(|z_{I^j}|)$ is contained in $z_{I^j}\neq0$ 
in particular, outside the ``critical" curve $\Gamma$ where superlogarithmic estimates hold by Theorem~\ref{t1.1} and Theorem~\ref{t1.1,5};
 thus $\di_b(\Pi_j\chi(|z_{I^j}|)$ are superlogaritmic multipliers.
\item[(b)] $\di_b\chi(|y_n|)\sim\cdot(h^{I^j}_{z_{I^j}})_j$ and hence it is by hypothesis a superlogarithmic multiplier.
\end{itemize}
Altogether, $\di_b\eta\contrazione$ are superlogarithmic multipliers. 
Remember that we are assuming that $(c_{ij}^h)$ are subelliptic multipliers.
Finally, supp$\,[\di_b,[\dib_b,\chi(z_{I^j})]]$ is contained in $z_{I^j}\neq0$ and $[\di_b,[\dib_b,\chi(y_n)]]\sim h^{I^j}_{z_{I^j},\overline{z_{I^j}}}$ are subelliptic multipliers; in conclusion, $[\di_b,[\dib_b,\eta]]$ are superlogarithmic multipliers. 
We can then apply Theorem~\ref{t1.3} and this completes the proof of Theorem~\ref{t1.4}

\hskip12cm $\Box$

\end{document}